\documentclass[conf]{new-aiaa}

\usepackage[english]{babel}
\usepackage{subcaption}
\usepackage{algorithm}

\usepackage{algpseudocode}
\usepackage{fancyhdr}
\usepackage{amsmath}
\usepackage{graphicx}
\usepackage{float}

\title{Rootfinding and Optimization Techniques for Solving Nonlinear Systems of Equations Arising from Cohesive Zone Models}
\author{Alberto Cattaneo\footnote{Graduate Student, Kahlert School of Computing.}}
\affil{Kahlert School of Computing, University of Utah, Salt Lake City, UT, 84112, United States of America}

\author{Varun Shankar\footnote{Assistant Professor, Kahlert School of Computing.}}
\affil{Kahlert School of Computing, University of Utah, Salt Lake City, UT, 84112, United States of America}

\author{M. Keith Ballard\footnote{Research Aerospace Engineer, Materials and Manufacturing Directorate, Composite Performance}}
\affil{Air Force Research Laboratory, Wright-Patterson AFB, OH, 45324, United States of America}

\begin{document}
\maketitle

\section{Introduction}

The ability to accurately predict the formation and evolution of cracks provides an important capability for the design, fabrication, and sustainment of composite materials. However, as researchers and engineers leverage increasingly complex composite material systems, modeling fracture becomes more difficult. A wide variety of approaches exist for modeling damage, but fundamental to all is finding solutions to highly nonlinear systems of equations.\\

Of the approaches for discretely modeling damage, this work focuses on the subset that use cohesive zone models. A cohesive zone model (CZM) is intended to capture the energy dissipated during crack opening while also directly accounting for a jump in the displacement field. A key advantage of using CZMs for cracks is that a stress singularity at the crack tip is avoided, which would introduce several challenges if considered. CZMs account for the energy dissipated over a characteristic length, called the ``fracture process zone''. There is a rich set of CZMs in the literature going back to the 1960s with the pioneering works by Dugdale\cite{CZMpaper}, Barrenblatt\cite{BARENBLATT196255}, and Hillerborg\cite{HILLERBORG1976773}. In early works, CZMs required identification of the fracture path \textit{a priori}, but the cohesive segment method was proposed in the early 2000's by Remmers et al. \cite{REMMERS03} to incrementally insert "cohesive segments" as cracks grew. This method did not require identification of a crack path \textit{a priori}, alleviating a critical limitation of earlier approaches. Since then, multiple methods have been proposed to account for the displacement jump along the evolving crack surfaces ranging from extended finite element methods (XFEM), which modify the basis functions to capture the displacement jump,\cite{HANSBO20043523,BALLARD2022114221,IARVE2011} to methods that automatically remesh the neighborhood of the crack\cite{CHEN2014104}.\\

While approaches to model the progression of fracture have received significant attention, methods to find the solution to the associated nonlinear equations have not. In general, nonlinear solution methods and optimization methods have a rich body of work spanning back to at least the first century\cite{Heath_2013}, providing the opportunity for advancement in the field of computational discrete damage modeling. In this paper, we explore the performance of established methods when applied to problems involving CZMs to identify promising methods for further improvement in this specialized application. We first use a simple 1D example problem with low degrees of freedom (DoF) to compare nonlinear solution methods, thereby allowing for both straightforward and intuitive visualization of the residual space and reasoning about the cause for each method's performance. We then explore the impact of higher DoF discretizations of the same problem on the performance of the solution methods. Finally, we discuss techniques to improve performance or to overcome limitations of the various methods.

\section{Finite Element Formulation}
In this section, we provide a brief overview of the governing equations for static equilibrium of a solid, deformable body and a corresponding finite element formulation for the 1D case. In the absence of body forces, the strong form of the governing equation is given by
\begin{align}\label{eq:strongform}
    \nabla \cdot \Tilde{\sigma} \left(x\right) = 0
\end{align}
where $\Tilde{\sigma}\left(x\right)$ represents the stress, generally a rank-2 tensor, at the point $x$. To facilitate a numerical solution, we introduce a test function  $\psi(\mathbf{x})$, integrate over the entire domain $\Omega$, and use integration by parts to obtain the weak form of the governing equation
\begin{align}\label{eq:weakform}
    \int{\Tilde{\sigma}(\mathbf{x})\cdot\nabla\psi(\mathbf{x})d\Omega} = \oint\psi(\mathbf{x})t_0(\mathbf{x})d\Gamma_0 + \oint \psi(\mathbf{x})t_c(\mathbf{x})d\Gamma_c
\end{align}
where $\Gamma_0$ is the external boundary of the domain $\Omega$, $t_0(\mathbf{x})$ is the traction applied along $\Gamma_0$, $\Gamma_c$ is the boundary introduced by cracks, and $t_c(\mathbf{x})$ is the traction that may exist along cracks due to the CZM used to model crack opening.\\

We construct the finite element formulation based on this weak form. We introduce a trial function to provide an approximation for the solution. As is commonly the case in finite element formulations, we use Lagrange polynomials compactly defined over each element. Therefore, the displacement, $u$, at any point is given by
\begin{align}\label{eq:interpolation}
    \mathbf{u}(\mathbf{x}) = \sum_i \theta_i(\mathbf{x})\mathbf{u}_i
\end{align}
where $\theta_i(\mathbf{x})$ defines the $i^{th}$ trial function and $\mathbf{u}_i$ represents the value of the displacement at the corresponding support point. Since the presented preliminary results are limited to 1D cases, we will restrict the present discussion of the formulation to 1D, though the theory is extensible to 3D. Additionally, while a test function can be selected from a different function space than the trial function, we use the same Lagrange polynomials for both the test and trial functions, thereby yielding a standard nodal (continuous) Galerkin finite element formulation. For convenience, we also assume the support points of the functions are enumerated such that $\theta_0$ and $\theta_N$ lie at $x_0$ and $x_N$, respectively. Finally, assuming linear elasticity and infinitesimal strain, Eqn.~\eqref{eq:weakform} can be expressed as
\begin{align}\label{eq:governing1}
    \int_{x_0}^{x_N} E(x)\frac{\partial\theta_i(x)}{\partial x}\frac{\partial\theta_j(x)}{\partial x}dx \: u_j = \theta_i(x)t_0(x)\biggr\rvert_{x_0}^{x_N} + \sum \theta_i(x)t_c(\Delta u)
\end{align}
where the domain extends from $x_0$ to $x_N$, $E(x)$ is the Young's modulus, and $t_c(\Delta u)$ is a function of the displacement jump $\Delta u$.\\

As introduced earlier, researchers have proposed a variety of CZMs, but for this study, we use a simple piece-wise linear model proposed in~\cite{turon2006} that defines the traction as
\begin{align}\label{eq:cohModel}
    t_c(\Delta u) = (1-d(\Delta u))K_p\Delta u
\end{align}
where $d(\Delta u)$ is a nonlinear damage parameter defined in~\cite{turon2006} and $K_p$ is a penalty parameter that should be sufficiently larger than $E(x)$.\\

The CZM introduces non-linear terms into the system of equations. Consequently, we define a residual function that measures the imbalance of forces between those experienced in the solid and those imposed at boundaries. For the presented preliminary results, we restrict "cracks" in our 1D model to locations at support points of $\theta_i$. Therefore, we describe a cohesive zone (CZ) by the pair of basis functions it connects. We define the residual as follows
\begin{align}\label{eq:resEqu}
    r_i = \int^{x_N}_{x_0} E(x)\frac{\partial\theta_i(x)}{\partial x}\frac{\partial\theta_j(x)}{\partial x}d\Omega \: u_j - \left( t_0(x_N)\delta_{iN} - t_0(x_0)\delta_{i0} \right) - \sum_{\{m,n\} \in C}(1-d(\Delta u ))K_p(u_m - u_n)(\delta_{im}-\delta_{in})
\end{align}
where $r_i$ is the residual at the $i^{th}$ support point, $\{m,n\}$ denotes the pair of basis functions that the respective CZ connects, and $C$ denotes the set of all CZs. The residual provides a measure of the imbalance of forces at each support point, and a valid solution exists when each term of the residual is 0.\\

Motivating this study, the terms in \eqref{eq:resEqu} can produce systems of equations with sharp nonlinearity and a poor condition number. The penalty stiffness for the CZM, $K_p$, introduces an artificial compliance, especially before the onset of damage. Consequently, $K_p$ is generally selected to be orders of magnitude larger than the stiffness of the surrounding material. Importantly, as the cohesive zone damage parameter $d(\Delta u)$ increases, the last term in \eqref{eq:resEqu} goes from a value that is very large to eventually zero. This behavior leads the Jacobian to change drastically when a change in the solution induces a significant change in the damage parameter. As a result, while CZMs are useful for predicting cohesive fracture, they can result in nonlinear systems of equations that are difficult to numerically solve.

\section{Methods}
We now describe our methods for solving the system of nonlinear equations corresponding to Eqn. \eqref{eq:resEqu}. Rewriting Eqn. \eqref{eq:resEqu} in a more general form that reflects its nature as a \emph{system} of equations, we obtain
\begin{align}\label{eq:main_solve}
    \mathbf{r}(\mathbf{u}) = K(\mathbf{u}) \: \mathbf{u} - \mathbf{b},
\end{align}
where $\mathbf{r}$ is the residual vector whose elements are $r_i$ from Eqn. \eqref{eq:resEqu}, $K: \mathbb{R}^n \to \mathbb{R}^{n \times n}$ is a nonlinear matrix-valued function, $\mathbf{u} \in \mathbb{R}^n$ is the vector of unknown displacements, and $\mathbf{b} \in \mathbb{R}^n$ contains known quantities corresponding to force terms. Our goal is to find $\mathbf{u}$ such that $\mathbf{r}(\mathbf{u}) = {\bf 0}$.\\

This section is divided as follows. In Section \ref{sec:list_methods}, we discuss the different nonlinear solution methods considered. We then present some common algorithms shared by all our algorithms in Section \ref{sec:common_alg}. Finally, we present our solvers in Section \ref{sec:solver_alg}.  In the remainder of this work, we occasionally shorten $\phi\left(\mathbf{r}(\mathbf{u}) \right)$ from \eqref{eq:phi} to $\phi(\mathbf{u})$ for clarity and brevity.

\subsection{Nonlinear solution methods considered}
\label{sec:list_methods}

Nonlinear solution methods are usually broken into two categories, optimization methods (OPT) and root-finding methods (RF). We consider eleven different methods that are optimization methods, root-finding methods, and hybrids of two.

\paragraph{Root-finding methods} We investigated root-finding methods based on the \emph{fixed-point iteration} $\mathbf{u}_{k+1} = \mathbf{g}\left(\mathbf{u}_{k}\right)$, where $g$ is some \emph{contractive function} and $\alpha$ is the step-size. Many popular numerical methods fall into this category including (but not limited to) the Picard method, the Newton method for rootfinding labeled Newton (RF) below, and the Broyden method. In our preliminary experiments, we restricted our attention to the Picard method, the Newton (RF) method, and two versions of Broyden's (so-called) ``Good method''~\cite{Broyden1965}. One version of the Broyden's Good Method, labeled Broyden below, approximates the Jacobian using rank-1 updates to some initial matrix (typically the identity matrix), while the other version, labeled Broyden (inv) below, approximates the inverse directly via the same strategy.  We briefly describe these below:

\begin{align}
    \text{Picard}: \ \mathbf{u}_{k+1} &= \mathbf{r}\left(\mathbf{u}_k\right), \label{eq:rm_picard} \\
    \text{Newton (RF)}: \ \mathbf{u}_{k+1} &= \mathbf{u}_k -  \left[\left.\mathbf{J}_{\mathbf{u}}\left(\mathbf{r} \right)\right|_{\mathbf{u} = \mathbf{u}_k}\right]^{-1} \mathbf{r}\left(\mathbf{u}_k\right), \label{eq:rm_newton} \\
    \text{Broyden}: \ \mathbf{u}_{k+1} &= \mathbf{u}_k -  \left[\left.\mathbf{J}'_{\mathbf{u}}\left(\mathbf{r} \right)\right|_{\mathbf{u} = \mathbf{u}_k}\right]^{-1} \mathbf{r}\left(\mathbf{u}_k\right), \label{eq:rm_Broyden} \\
    \text{Broyden (inv)}: \ \mathbf{u}_{k+1} &= \mathbf{u}_k -  \left[\left.\mathbf{J}_{\mathbf{u}}\left(\mathbf{r} \right)\right|_{\mathbf{u} = \mathbf{u}_k}\right]^{-1'} \mathbf{r}\left(\mathbf{u}_k\right), \label{eq:rm_BroydensInverse}
\end{align}

where $J_{\mathbf{u}}(.)$ is the Jacobian of the argument with respect to the vector $\mathbf{u}$; in this case, the Jacobian is a matrix-valued function $J: \mathbb{R}^n \to \mathbb{R}^{n\times n}$. $J'$ and $J^{-1'}$ are estimates of the Jacobian, and the inverse of the Jacobian, respectively. The convergence of these methods hinges on the function $\mathbf{g}$ being contractive. Unfortunately, for problems involving CZMs, $\mathbf{g}$ can become non-contractive, causing convergence to stall.

\paragraph{Optimization methods} Optimization methods solve Eqn. \eqref{eq:main_solve} by transforming it into an optimization problem whose minimum matches the solution of Eqn. \eqref{eq:main_solve}. The optimization methods considered in this work use a scalar-valued \emph{objective function} $\phi\left(\mathbf{r}(\mathbf{u})\right)$ given by
\begin{align}\label{eq:phi}
\phi\left(\mathbf{r}\left(\mathbf{u}\right) \right) = \frac{1}{2}\left\|\mathbf{r}\left(\mathbf{u}\right) \right\|_2^2.
\end{align}
The optimization methods we considered require both the Jacobian and Hessian of $\phi$ with respect to $\mathbf{u}$. The Jacobian $\mathbf{J}_{\mathbf{u}}(\phi)$ can be obtained by an application of the chain rule for function composition:
\begin{align}\label{eq:jacob_opt}
    \mathbf{J}_{\mathbf{u}}\left(\phi\left(\mathbf{r}(\mathbf{u})\right) \right) = \left[\mathbf{J}_{\mathbf{u}} \left(\mathbf{r}\right)\right]^T \nabla_{\mathbf{r}} \phi = \left[\mathbf{J}_{\mathbf{u}} \left(\mathbf{r}\right)\right]^T \mathbf{r},
\end{align}
where $\nabla_{\mathbf{r}}$ is the gradient with respect to $\mathbf{r}$. The Hessian of $\phi$, denoted by $H_{\mathbf{u}}(\phi)$ is given by two applications of the chain rule. However, for the tests considered herein, we approximate the Hessian as
\begin{align}\label{eq:hess_opt_approx}    \mathbf{H}_{\mathbf{u}}\left(\phi\left(\mathbf{r}(\mathbf{u})\right) \right) \approx \left[\mathbf{J}_{\mathbf{u}} (\mathbf{r})\right]^T \left[\mathbf{J}_{\mathbf{u}} (\mathbf{r})\right],
\end{align}
which effectively drops a term involving second derivatives from the exact Hessian; this approximation allows for efficient computation of the approximate Hessian from the Jacobian. In this work, we considered three quasi-Newton optimization methods, all of which involve using approximations to the Hessian, either given by Eqn. \eqref{eq:hess_opt_approx} or via some other technique. In this group are methods such as the ADAM \cite{Duchi2011,Diederik2014} family of methods, L-BFGS\cite{Matthies1979} and BFGS all described below. They are as such:
\begin{align}
    \text{ADAM}: \ \mathbf{u}_{k+1} &= \mathbf{u}_k -  \left[\left.aggregate(\mathbf{\nabla_{\mathbf{r}}}\left(\mathbf{u} \right)\right|_{\mathbf{u} = \mathbf{u}_k})\right]\alpha
    \label{eq:opt_adam} \\
    \text{ADAGRAD}: \ \mathbf{u}_{k+1} &= \mathbf{u}_k -  \left[\left.aggregate(\mathbf{\nabla_{\mathbf{r}}}\left(\mathbf{u} \right)\right|_{\mathbf{u} = \mathbf{u}_k}^2)\nabla_{\mathbf{r}}\right]\alpha, \label{eq:opt_adagrad} \\
    \text{BFGS}: \ \mathbf{u}_{k+1} &= \mathbf{u}_k - H_u'(\phi)\nabla_{\mathbf{r}}, \label{eq:opt_bfgs} \\
    \text{L-BFGS}: \ \mathbf{u}_{k+1} &= \mathbf{u}_k - H_u''(\phi)\nabla_{\mathbf{r}}, \label{eq:opt_lbfgs}
\end{align}
For L-BFGS and BFGS we use $H'$ and $H''$ to denote two different ways of estimating the Hessian using either the entire history or a finite number of past iterations.

\paragraph{Hybrid Methods}
We also tested two trust-region quasi-Newton methods, both of which maintain a constraint on the size of the optimization step. The first was a modification of Powell's Dogleg method~\cite{PowelDogleg}, in which we linearly combined a gradient descent step on $\phi$ using Eqn. \eqref{eq:jacob_opt} with a rootfinding Newton step using Eqn. \eqref{eq:rm_newton}. Our modified Dogleg method uses a backtracking line search (Algorithm \ref{alg:ls}) to determine the size of the gradient descent step and is given in Algorithm \ref{alg:dl}.
It assumes that the Newton step, if within the trust region, is likely to be a better estimate. If the Newton step is outside the trust region, we form a linear combination of the Newton and gradient descent step that is guaranteed to stay within the trust region. As is typical for trust-region methods, we used a local quadratic model to determine whether a step in a given direction is acceptable and by extension whether the trust-region should be contracted or expanded. The model itself given a direction $\mathbf{p}$ is as follows
\begin{align}\label{eq:quadModel}
m(\mathbf{u})=\phi\left(\mathbf{r}\left(\mathbf{u}\right)\right) +\nabla_{\mathbf{u}} \phi\left(\mathbf{r}\left(\mathbf{u}\right)\right)^T\mathbf{p} + \frac{1}{2}\mathbf{p}^T H_{\mathbf{u}}\left(\phi\left(\mathbf{r}(\mathbf{u})\right)\right) \mathbf{p}
\end{align}
The second trust-region method that we tested was the Steihaug method~\cite{SteihaugMethod}, which is typically viewed as a modified conjugate gradient method that ensures steps stay within a trust region. Our modification to the Steihaug method was dicated by practicality. In situations where Steihaug could not find a suitable step size within a set number of sub-iterations, we simply reset the trust region and returned the current step size. The modified Steihaug method is given in Algorithm \ref{alg:st}. In both cases, the trust region radius $\delta$ is a hyperparameter that is set by the user. However, in Algorithm \ref{alg:at}, we present a simple technique to modify $\delta$, allowing even greater flexibility. \\

\textbf{Note}: While the methods can be divided as discussed above, another important distinction between methods is whether the Jacobians (and/or Hessians) are exactly calculated or estimated/approximated. Newton (RF), Steihaug, and Dogleg use the exact Jacobian. In contrast, Broyden (both versions), (L-)BFGS, ADAM, and ADAGRAD all use approximate gradients, Jacobians, or Hessians. These methods use updates based on the steps taken in order to create either estimates of first and second order information (Broyden, BFGS) or modifications of the gradient based on previous steps (ADAM, ADAGRAD). This alternate classification proves to be important in Section \ref{sec:scalingTrends}.

\subsection{Common algorithms}
\label{sec:common_alg}
%
%
We now present two useful algorithms for use within our solvers: a backtracking line search (Algorithm \ref{alg:ls}) and an algorithm for modifying the trust region radius $\delta$ (Algorithm \ref{alg:at}).
\begin{algorithm}[!htpb]
\caption{Backtracking line search}\label{alg:ls}
\begin{algorithmic}
\Require{$\phi$ is differentiable}
\Require{$\nabla \phi$ is known}
\Require{Current position $\mathbf{u}$, search direction $\mathbf{p}$, and maximum step size $\alpha_0$}
\Require{Control parameters $\tau$ and $c$}
\State $q = \nabla_{\mathbf{u}} \phi(\mathbf{u})^T \mathbf{p} $
\State $t = -c q$
\State $j = 0$
\While{$\phi(\mathbf{u}) - \phi(\mathbf{u} + \alpha_j\mathbf{p}) < \alpha_j t$ and Wolfe conditions not met}
\State $\alpha_{j+1} = \alpha_{j}\tau$
\State $j = j+1$
\EndWhile
\end{algorithmic}
\end{algorithm}
The Wolfe conditions are a set of inequalities that aid inexact line searches by constraining the step size $\alpha$ \cite{armijo,wolfe1971}. In Algorithm \ref{alg:ls}, we use two Wolfe conditions, the Armijo condition and the curvature condition. A step size meets the Armijo condition if it satisfies
\[
\phi(\mathbf{u}_k + \alpha\mathbf{p}_k) \leq \phi(\mathbf{u}_k) + c_1\alpha_k\mathbf{p}_k^T\nabla \phi(\mathbf{u}_k).
\]
It is said to satisfy the curvature condition if it meets the following equality:
\[
-\mathbf{p}_k^T\nabla \phi(\mathbf{u}_k + \alpha\mathbf{p}_k) \leq -c_2\mathbf{p}_k^T\nabla \phi(\mathbf{u}_k).
\]
It is also important to note that other objective functions can be used in the line search, not just $\phi$. For instance, in Algorithm \ref{alg:st} we use the model $m(\mathbf{u})$ from Eqn. \eqref{eq:quadModel}. 
\begin{algorithm}[H]
\caption{Adjusting $\delta$}\label{alg:at}
\begin{algorithmic}
\State $j= 0$

\State $p_{\phi} = \phi(\mathbf{u}_{k-1}) - \phi(\mathbf{u}_k) $
\State $m_k = m(u_{k+1}-u_k)$ \eqref{eq:quadModel}
\State $p_{m} = \phi(\mathbf{u}_{k-1}) - m_k $
\State $p = \frac{p_{\phi}}{p_{m}}$
\State  $j= j + 1$
\If{$p < 0.25$}
\State $\delta_{k+1} = 0.25\delta_{k}$
\EndIf
\If{$p > 0.75$}
\State $\delta_{k+1} = \min\left(2\delta_{k} ,\delta_{max}\right)$
\EndIf
\If{ $j= 5$}
\State $\delta_k = \delta_0$
\State $j = 0$
\EndIf
\end{algorithmic}
\end{algorithm}
After every update to the guess $\mathbf{u}$, we also update the trust region. This is described at the end of each of the solver algorithms as "Adjust $\delta_k$ for $\delta_{k+1}$". We found in our experiments that it was beneficial to periodically reset $\delta$. We employed this approach in order to prevent the trust region from shrinking to the empty set (which could happen if $\phi$ does not match its local quadratic model sufficiently). These details are given in Algorithm \ref{alg:at}.

\subsection{Solver algorithms}
\label{sec:solver_alg}
We now present the algorithm details for our modified Powell's dogleg method (Algorithm \ref{alg:dl}) and our modified Steihaug method (Algorithm \ref{alg:st}). These run until convergence is reached as defined by the norm of the residual decreasing to some tolerance $\epsilon$. 

It is also worth remarking on detecting divergence or stalling in any of these algorithms. 
Because of the existence of necessary sub-optimal-steps for several methods, checking for an increase in the residual is insufficient to check for divergence or stalling. Unlike checking for a simple increase in residual, a way of checking for non-convergence involves keeping track of the history of the residual norms, finding the curve through the residual norms using a polynomial, and checking whether that curve has a negative or positive trend. However, this procedure does require the introduction of another parameter specifying the number of previous residual norms kept. Additionally, some methods, such as Broyden, will also fail if their convergence stalls near the root/minimum; this situation can be prevented by checking if a denominator in the algorithm approaches 0 and allowing for early termination if this occurs.
\begin{algorithm}[H]
\caption{Modified Powell's dogleg}\label{alg:dl}
\begin{algorithmic}
\Require{Initial guess $\mathbf{u}_0$, initial trust region radius $\delta_0$, maximum trust region radius $\delta_{max}$}
\Require{The Jacobian $\mathbf{J}_k = \left.\mathbf{J}_{\mathbf{u}} \left(\mathbf{r}(\mathbf{u})\right)\right|_{\mathbf{u} = \mathbf{u}_k}$ is known}
\Require{The minima of $\phi$ is equal to a root of $\mathbf{r}$}
\State $k = 0$

\While{$\|\mathbf{r}(\mathbf{u}_k)\|_2 > \epsilon$}
\State $\mathbf{p}_g = \left.\nabla_{\mathbf{u}} \phi(\mathbf{u})\right|_{u=u_k} = -\mathbf{J}_k^T \mathbf{r}(\mathbf{u}_k)$
\State $\alpha_k =$ Backtracking Line Search ($\mathbf{p}_g$,$\mathbf{u_k}$) (\ref{alg:ls})

\If{$\|\alpha \mathbf{p}_g\|_2 > \delta_k$} \Comment{(OPT) Gradient descent step}
\State $\mathbf{u}_{k+1} = \mathbf{u}_k + \beta \mathbf{p}_g$ with $\beta$ for which $\|\beta \mathbf{p}_g\|_2 = \delta_k$
\EndIf

\State $\mathbf{p}_n = \mathbf{J}^{-1} \mathbf{r}(\mathbf{u}_k)$ \Comment{(RF) Newton step}
\If{$\|\mathbf{p}_n\|_2 \leq \delta_k$}
\State $\mathbf{u}_{k+1} = \mathbf{u}_k + \mathbf{p}_n$
\EndIf
\State $a = \|\mathbf{p}_n - \mathbf{p}_g\|_2^2$
\State $b = \mathbf{p}_g^T (\mathbf{p}_n - \mathbf{p}_g)  $
\State $c = \|\mathbf{p}_n\|_2^2 - \delta_k^2$

\State \(  \tau = \frac{-b \pm \sqrt{b^2 - 4ac}}{2a}   \)
\State $\mathbf{p}_m = \mathbf{p}_g - \tau(\mathbf{p}_n-\mathbf{p}_g)$\Comment{Cauchy point adjustment}

\If{$\|\mathbf{p}_m\|_2 \leq \delta_k$} 
\State $\mathbf{u}_{k+1} = \mathbf{u}_k + \mathbf{p}_m$
\Else{ $\mathbf{u}_{k+1} = \mathbf{u}_k + \beta \mathbf{p}_m$ with $\beta$ for which $\|\beta \mathbf{p}_m \|_2 = \delta_k$}
\EndIf

\State $k=k+1$
\State Adjust $\delta_k$ for $\delta_{k+1}$ (\ref{alg:at})
\EndWhile
\end{algorithmic}
\end{algorithm}

\begin{algorithm}[!b]
\caption{Modified Steihaug}\label{alg:st}
\begin{algorithmic}[1]
\Require{Initial guess $\mathbf{u}_0$, initial trust region radius $\delta_0$, maximum trust region radius $\delta_{max}$ and the number of iterations to reset the trust region}
\Require{The Jacobian $J_k = \left.J_{\mathbf{u}} \left(\mathbf{r}(\mathbf{u})\right)\right|_{\mathbf{u} = \mathbf{u}_k}$ and the approximate Hessian $H_k = \left.H_{\mathbf{u}}(\phi(\mathbf{u}))\right|_{\mathbf{u} = \mathbf{u}_k}$ are known}
\Require{The minima of $\phi$ is equal to a root of $\mathbf{r}$}

\State $n = len(\mathbf{u}_0)$
\State $k = 0$
\algstore{bkbreak}
\end{algorithmic}
\end{algorithm}

\begin{algorithm}[!h]
\begin{algorithmic}
\algrestore{bkbreak}
\While{$\|\mathbf{r}(\mathbf{u}_k)\|_2 > \epsilon$}
\State $j = 0$ 
\State $\mathbf{d}_k = \nabla_u \phi(\mathbf{u}_k)$
\State $\mathbf{r}_k = -\mathbf{d}_k$
\While{$j < n-1$}
\State $c_j = \mathbf{d}_k^T H_k^T\mathbf{d}_k$
\If{$c_j \leq 0$}
\State $\beta\mathbf{d}_k$ with $\beta$ such that $\|\beta\mathbf{d}_k\|_2 = \delta _k$
\State $\mathbf{u}_{k+1} = \mathbf{u}_k + \beta\mathbf{d}_k$
\State Exit $j$ loop
\EndIf

\State $\alpha_k =$ Backtracking Line Search ($\mathbf{d}_k$,$\mathbf{u_k}$) (\ref{alg:ls})
\State $\mathbf{s}_j = \alpha_k\mathbf{d}_j$
\If{$\|\mathbf{s}_j\|_2 \geq \delta_k$}
\State $\beta\mathbf{s}_k$ with $\beta$ such that $\|\beta\mathbf{s}_k\|_2 = \delta _k$
\State $\mathbf{u}_{k+1} = \mathbf{u}_k + \beta\mathbf{s}_k$
\State Exit $j$ loop
\EndIf

\State $\mathbf{r}_k^{new} = \mathbf{r}_k + \alpha_j H_k \mathbf{d}_k$

\If{$\|\mathbf{s}_j\|_2 \leq \delta_k$}
\State $\mathbf{u}_{k+1} = \mathbf{u}_k + \mathbf{s}_k$
\State Exit $j$ loop
\EndIf

\State $\gamma_j = \frac{ \left(\mathbf{r}^{new}_k\right)^T \mathbf{r}^{new}_k}{(\mathbf{r}_k)^T \mathbf{r}_k}$
\State $\mathbf{d}_{k} = \mathbf{r}^{new}_k + \gamma\mathbf{d}_k$
\State $\mathbf{r}_{k} = \mathbf{r}^{new}_k$

\State $j = j+1$
\EndWhile

\State Adjust $\delta_k$ for $\delta_{k+1}$ (\ref{alg:at})
\State $\mathbf{u}_{k+1} = \mathbf{u}_k + \mathbf{s}_k$ \Comment{Optionally, if no step is found in the $j$ loop}
\State $k = k+1$
\EndWhile
\end{algorithmic}
\end{algorithm}
Because the time required for a single iteration in a method may differ, Table \ref{tab:allCosts} lists estimates of the cost of 1 iteration of each of the methods discussed. This table does not include the cost of doing the line-search. The O(n) complexity notation provides some information regarding the time required in the worst case, but each method will be sensitive to well-optimized implementations. Furthermore, where applicable, we use GMRES with no restarts and an internal tolerance of $10^{-5}$ as the linear solver in lieu of computing the inverse directly.

\begin{table}[H]
\caption{\label{tab:allCosts} Cost of 1 iteration for each of the methods tested. }
\centering
\begin{tabular}{l|r|r|r|r}
Method & Residuals & Jacobians & Solves & O(n) \\\hline
Dogleg* & 6 & 6 & 0 & O($n^2$)\\  
Steihaug* & 6 & 6 & 0 & O($n^2$)\\
Newton (RF) & 1 & 1 & 1 & O($n^3$)\\  
L-BFGS-B* & 1 & 3  & 0 & O($n^2$)\\
Broyden & 1 & 0  & 1 & O($n^3$)\\
Broyden (inv) & 1 & 0  & 0 & O($n^2$)\\
ADAM* & 3 & 4  & 0 & O($n^2$)\\
ADAGRAD* & 3 & 4 & 0 & O($n^2$)
\end{tabular}
\end{table}

\section{Results}
\label{sec:results}
We tested each of the nonlinear methods considered for two different boundary value problems (BVPs). For the first BVP, the equilibrium solution requires the cohesive zone to go from the initial undamaged condition to completely open, sustaining no traction. We hence refer to this case as ItC (initial to complete) in the following diagrams. For the second BVP, an equilibrium solution exists such that the damage parameter of the cohesive zone is between 0 and 1, corresponding to a partially damaged state. We refer to this case as ItP (initial to partial). As mentioned, we focused on 1D problems in order to explore characteristic boundary value problems while gaining an full understanding of the behavior exhibited by the different solvers. Importantly, the two cases each highlight an important class of BVPs. ItP is designed to be well-behaved with regard to fixed-point methods (Newton and Picard) while ItC is designed to fail with fixed-point methods. In our experiments, we defined convergence as the 2-norm of the residual having decreased to less than $10^{-6}$.

\subsection{Coarse Elasticity Problem with a Single Cohesive Zone}
\begin{table}[!htpb]
\caption{\label{tab:allresults} Number of iterations to convergence on the ItP and ItC cases with a coarse mesh. A value of $\infty$ indicates non-convergence.}
\centering
\begin{tabular}{l|r|r}
Method & Iterations (ItP) & Iterations (ItC) \\\hline
Dogleg & 2 & 83\\
Steihaug & 38 & 83 \\
Newton & 2 & $\infty$\\
L-BFGS-B & 12 & 16\\
ADAM & 719 & 1206 \\
ADAGRAD & 172 & 352\\
Picard & 2 & $\infty$\\
Broyden & 11 & 11\\
Broyden (inv) & 11 & 11\\
BFGS & $\infty$ & $\infty$
\end{tabular}
\end{table}
\begin{figure}[!htpb]
    \centering        
    \includegraphics[height=2.3in]{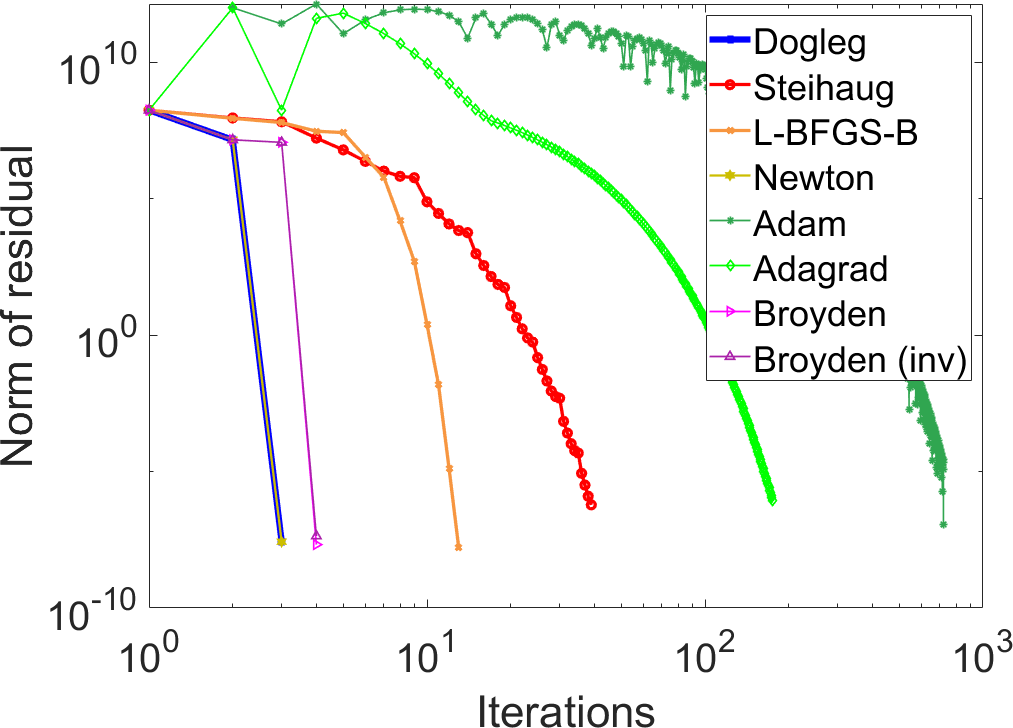}
    \includegraphics[height=2.3in]{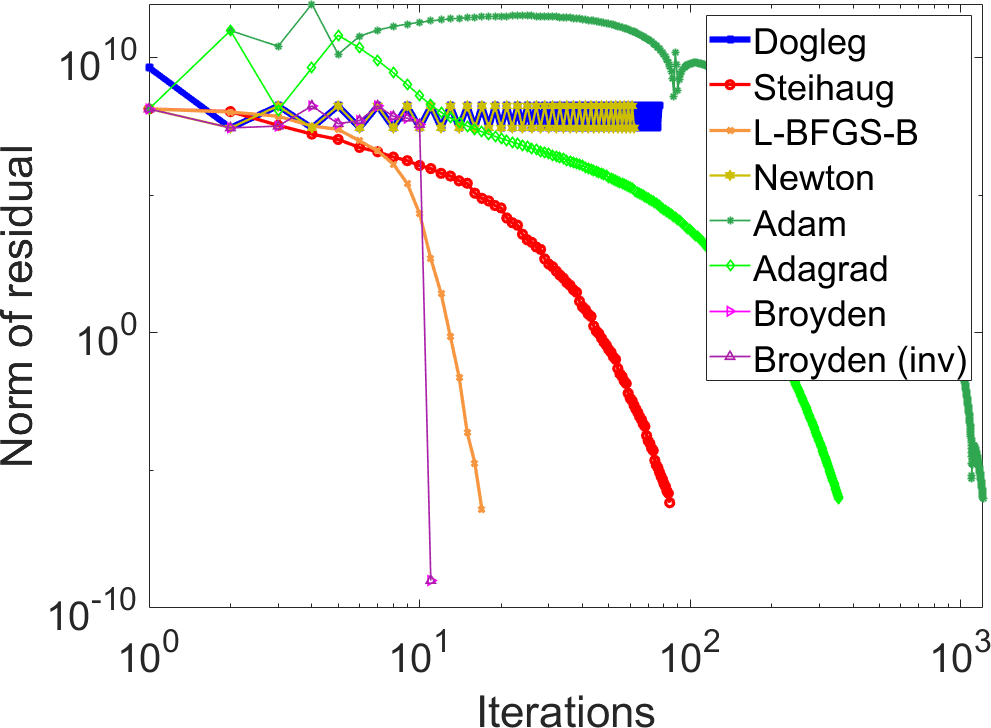}
    \caption{Residual norms as a function of iteration number for case ItP (left) and case ItC (right).} 
    \label{fig:res_norms}
\end{figure}
In order to get a clear picture of the nature of the problem, we began with the simplest version of the problem involving two linear elastic elements with a single cohesive zone element between them. The results are summarized in Figure \ref{fig:res_norms} and Table \ref{tab:allresults}. Table \ref{tab:allresults} shows the number of iterations each method required for each of the two BVPs. While all methods converged for the ItP case, the Newton, Picard, and Dogleg methods did so in two iterations. The Dogleg method  was able to use two Newton steps for this case, not requiring any gradient descent steps. L-BFGS converged faster than the trust-region Newton methods but slower than Newton, Picard, or Dogleg. For the ItC case, Newton and Picard failed to converge, while the trust-region methods (Dogleg and Steihaug) both converged in 83 iterations. L-BFGS appeared to perform best on the ItC case. We believe this to be due to a combination of factors including its nature as a quasi-Newton method (faster convergence than gradient-based methods), built-in regularization due to the use of approximate derivatives, and stabilization through the use of a line search. ADAM and ADAGRAD did converge on both problems, but took a very large number of iterations.

\begin{figure}[!htpb]
    \centering    
    \includegraphics[height=2.3in]{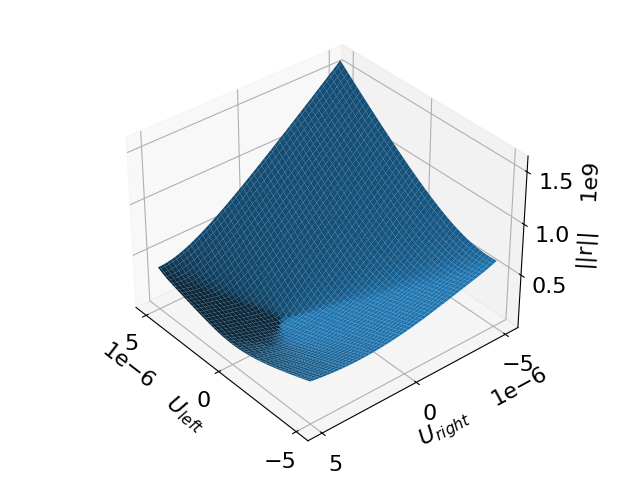}
    \includegraphics[height=2.3in]{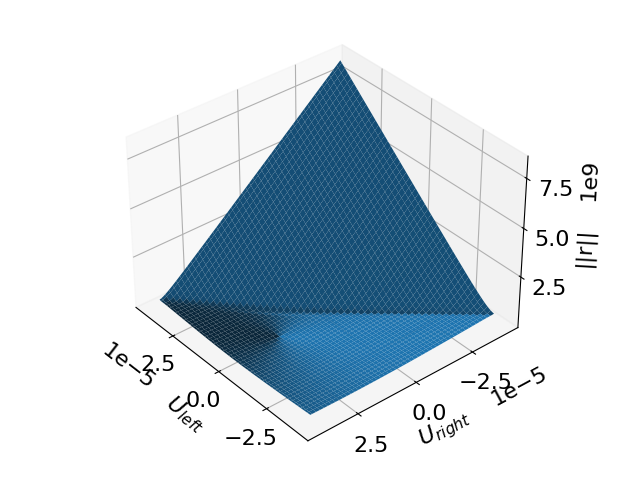}
    \caption{Residual surfaces for case ItP (left) and case ItC (right).}
    \label{fig:res_surf}
\end{figure}

Figure \ref{fig:res_norms} illustrates the results from Table \ref{tab:allresults} more starkly. The results for case ItP in Figure \ref{fig:res_norms} (left) show the rapid convergence of the fixed-point methods and the Dogleg method (which reverted to Newton); nevertheless, all methods converge, with L-BFGS offering performance closer to the fixed-point methods than to the trust-region ones. On the other hand, Figure \ref{fig:res_norms} (right) shows that the residual norms for Picard and Newton oscillated without ever converging, trapped in the trough in Figure \ref{fig:res_surf}. In contrast, the trust-region methods worked their way out of the trough. In fact, we found that they spent most of their iterations near the root (not shown).
\begin{figure}[!tpb]
    \centering        \includegraphics[width=0.5\linewidth]{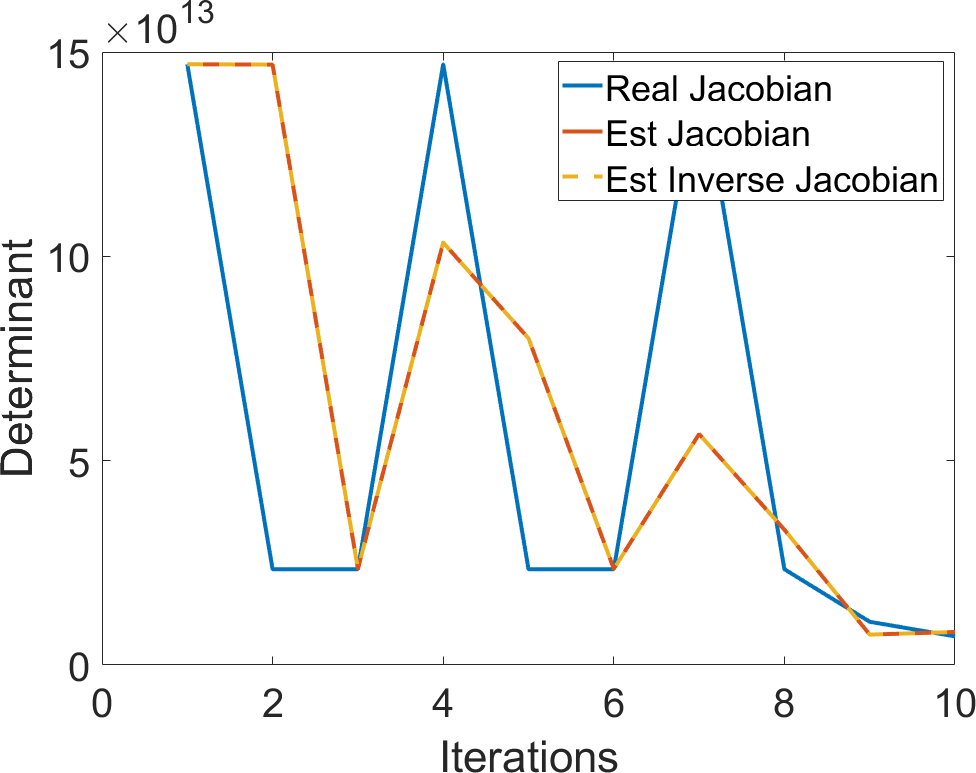}
    \caption{Estimated and exact Jacobian condition number during iterations with Broyden's method for ItC case.}
    \label{fig:condCompare}
\end{figure}
In general, the trend is that methods that rely on finding the minimum of the objective function perform more consistently, with L-BFGS doing especially well. We suspect this is because the transformed version of the function that we are optimizing over is better behaved than the original function, and therefore, the optimization methods are less likely to get stuck. Both variants of Broyden also performed well despite being a root-finding method, which may be attributed to its continual low-rank updates to the approximated Jacobian. To illustrate this, Figure \ref{fig:condCompare} shows the condition number of the estimated and exact Jacobian (or their inverses) during Broyden iterations. Importantly, changes in the condition number of the estimated Jacobian appear to lag behind those in the exact Jacobian and have decreasing magnitude over iteration count; this could potentially help avoid oscillations that can lead to non-convergence.\\

Several methods showcase behavior in which a step is taken that causes the residual to increase before the next iteration results in a decrease. These "sub-optimal steps" can provide a mechanism by which difficult residual surface geometry can be circumvented. Because methods like ADAM and L-BFGS keep a history of steps, sub-optimal steps can provide useful information to form robust estimates of the Jacobian. Furthermore, methods that create estimates of first and second order information \textit{must} take poor steps early in order to create the estimates which are then improved. The methods all eventually may decrease the residual after sufficient iterations, despite an early appearance of divergence.\\
Figure \ref{fig:res_surf} shows residual surfaces for both cases, where $u1$ and $u2$ refers to the displacement of the node on the left and right side of the cohesive zone, respectively. The residual for case ItP is simple and well-approximated by a quadratic. On the other hand, the residual for case ItC exhibits more complicated topology: two distinct regions separated by a ``trough'' which needs to be crossed to reach the root. This trough corresponds to an attractor for fixed point methods, which in our case corresponded to oscillatory behavior without any convergence. However, the residual surface for case ItP does not exhibit any complications, and indeed, the results in Table \ref{tab:allresults} our results show that fixed-point methods consequently perform very well. This illustrates the complexity present even for the simplest problem involving a cohesive zone model.

\subsection{Fine Elasticity Problem with a Single Cohesive Zone}
\label{sec:scalingTrends}
In this section, we present our experiments with using a finer mesh for the elasticity problem to explore the behavior of the methods as the linear portion of the system of equations increases. Instead of two linear elastic elements, we used a mesh consisting of 64 linear elastic elements with a single cohesive zone in the same location as the previous study. Importantly, methods that involve a linear solve with the exact Jacobian, such as Newton (RF) method, should exhibit similar behavior as observed for the coarse problem. However, methods that use estimates or avoid linear solves are anticipated to behave differently as more linear elements are introduced. \\
\begin{table}[!htpb]
\caption{\label{tab:allresults_fine} Number of iterations to convergence on the ItP and ItC cases with a finer mesh. A value of $\infty$ shows non-convergence due to divergence, while * indicates convergence but not to the required precision.}
\centering
\begin{tabular}{l|r|r}
Method & Iterations (ItP) & Iterations (ItC) \\\hline
Dogleg & 7* & >100\\
Steihaug & >100 & >100 \\
Newton & 6* & $\infty$\\
L-BFGS & >100 & >100\\
ADAM & >100 & >100 \\
ADAGRAD & >100 & >100\\
Picard  &$\infty$ & $\infty$ \\
Broyden & 18* & 31\\
Broyden (inv) & 15* & 30*\\
BFGS & $\infty$ & $\infty$
\end{tabular}
\end{table}
Table \ref{tab:allresults_fine} shows the number of iterations each method required for the two BVPs. Unlike the results from the previous section for the coarse mesh, no single family of methods demonstrated a clear advantage. The increase in degrees of freedom due to a finer mesh generally led to substantial increase of the number of iterations required. Furthermore, L-BFGS, ADAM, and ADAGRAD exhibit impractically slow convergence. The best performing methods are both versions of the Broyden method. Overall, it appeared that performance on this finer mesh with more DoFs related to whether the method used exact or estimated Jacobians/Hessians, rather than whether the method relied on rootfinding or optimization. \\

\begin{figure}[!htpb]
    \centering
    \includegraphics[height=2.3in]{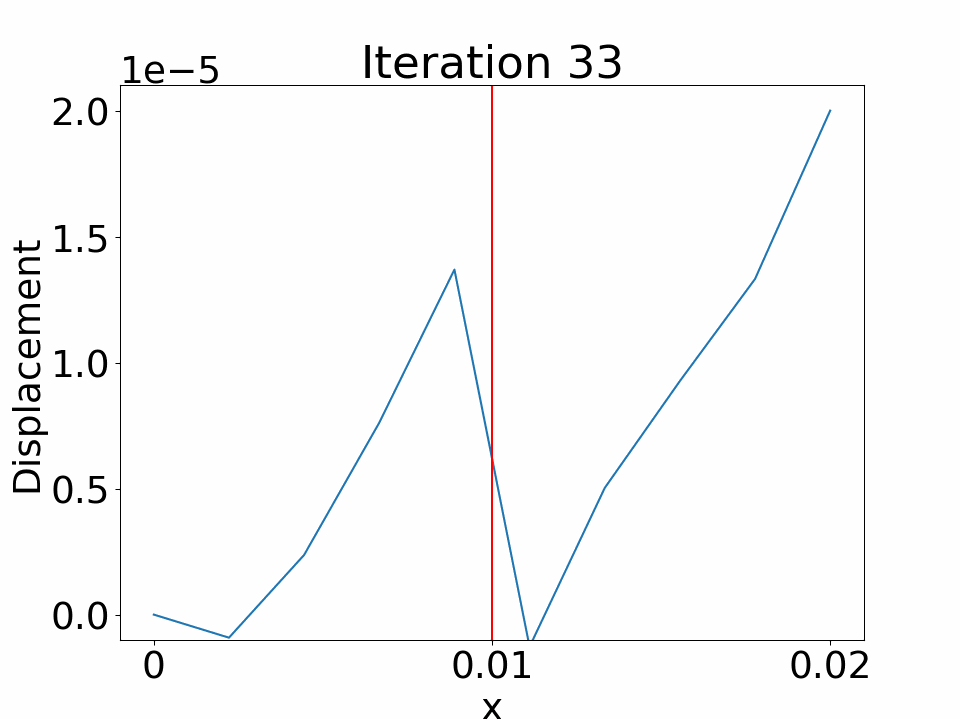}
    \includegraphics[height=2.3in]{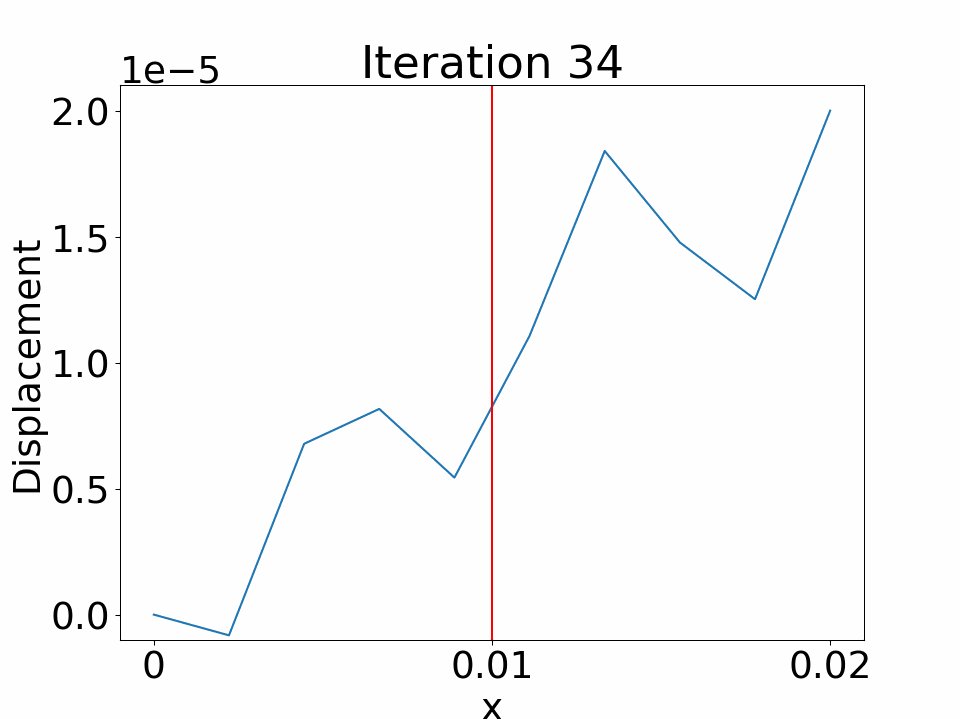}
    \caption{Example of instability across iterations of the L-BFGS-B method for the ItC case. The red line is the location of the cohesive zone.} 
    \label{fig:case1LBFSGinstability}
\end{figure}

All methods exhibited one of two problems: a loss of precision or divergence. Both Newton (RF) and Broyden stalled as they approached the root but before a sufficiently low residual norm was reached. Newton (RF), ADAGRAD, and ADAM stalled as a consequence of the aforementioned ill-conditioning of the exact Jacobians at each iteration, while Broyden stalled because the method is unable to find all of the spectral information and thus cannot get any closer. The methods that diverge do so due to numerical instability, an example of which can be seen in Figure  \ref{fig:case1LBFSGinstability}. As was mentioned before, several methods require sub-optimal steps to be taken. During these steps, L-BFGS can encounter overflows due to the magnitude of the steps it takes. On the other hand, the Dogleg and Steihaug methods move in a poor direction due to numerical instability, leading to divergence.\\

\begin{figure}[!htpb]
    \centering  
    \includegraphics[height=2.3in]{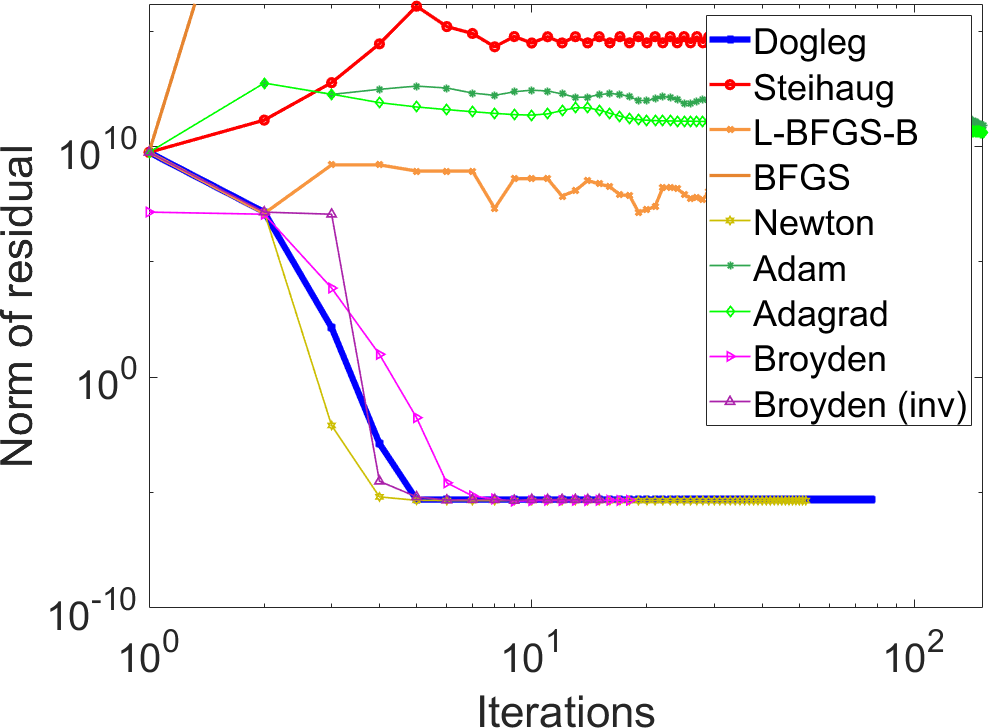}
    \includegraphics[height=2.3in]{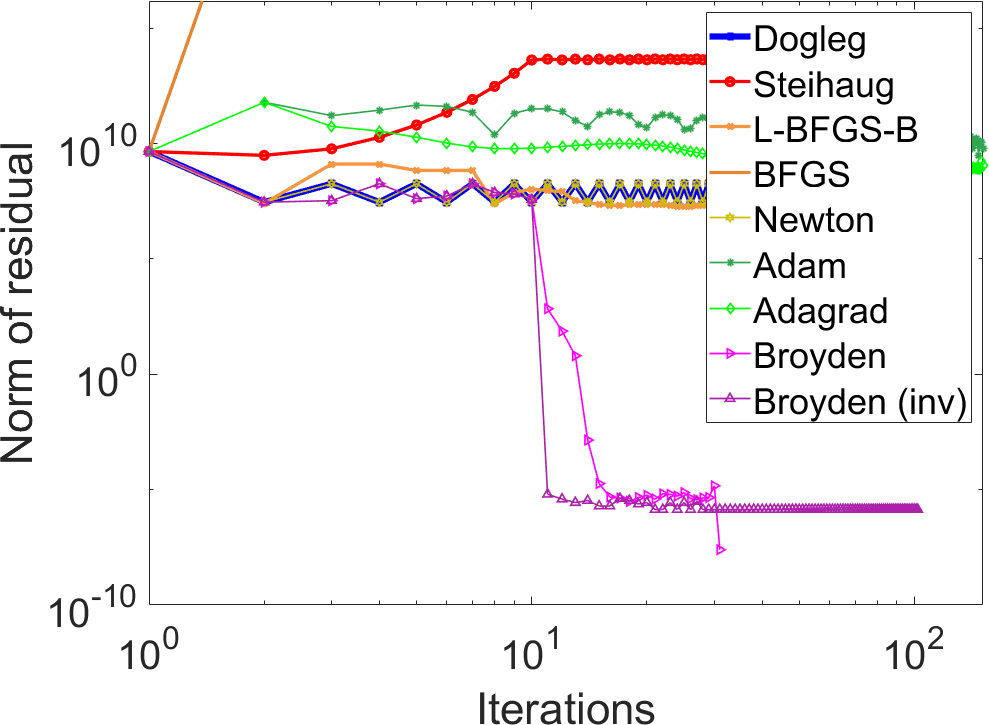}
    \caption{Residual norms as a function of iteration number for case ItP (left) and case ItC (right) with a finer mesh.}
    \label{fig:res_norms_64}
\end{figure}

Figure \ref{fig:res_norms_64} shows the history of the residual norm throughout iteration for each method. Newton (RF) exhibited performance similar to that shown for the coarse mesh. For the ItP case, the residual decreased rapidly, but unlike the previous results for the coarse mesh, it stalled near the solution. We expect that this behavior is due to the linear solver selected, and if an exact linear solve was used, convergence would more closely resemble the results for the coarser mesh. However, for the ItC case, Newton (RF) failed in similar manner to the coarse mesh results. Unlike for the coarse mesh, the Steihaug method exhibited an early increase in the residual followed by plateau. It may converge eventually, but the number of iterations would be impractical. Dogleg demonstrates similar performance to Newton (RF), since the Newton steps are unable to reach the solution and the non-Newton steps do not take it in the correct direction. The Broyden methods perform consistently best for both speed and precision. While the two variants performed almost identically for the coarse mesh, differences appear for the finer mesh with Broyden (inv) performing slightly better in the neighborhood of the solution. L-BFGS performed much worse with a finer mesh, and BFGS diverged quickly. Both ADAM and ADAGRAD do not diverge, but they took an impractical number of iterations to reduce the residual. We have excluded Picard in this case because it becomes completely unstable and does not fit in the plot.

\subsection{Conclusions}
In this paper, we explored how a variety of nonlinear solution methods performed for two types of 1D BVPs involving a single cohesive zone (CZ). The first case involves boundary conditions designed to have a solution where the CZ has a damage variable between 0 and 1 (referred to as the initial-to-partial or ItP case), while the second case has a different set of boundary conditions such that the only solution corresponds to the CZ being completely open (referred to as the initial-to-complete or ItC case). While the BVPs are far less complex than typical engineering problems, the systems of equations exhibit strong nonlinearity because the cohesive zone has a significant effect on the total solution, making the BVPs useful for a comparison of methods. The nonlinear solution methods considered included root-finding methods (Picard, Newton, and two variants of Broyden), optimization methods (ADAM, ADAGRAD, BFGS, and L-BFGS), and hybrid methods (Dogleg and Steihaug). \\

When using a coarse mesh consisting of two linear elastic elements connected by a cohesive element, the performance of methods was generally divided along rootfinding (RF) and optimization (OPT) methods. RF quickly converged for the ItP case but were unable to converge for the ItC cases with the exception of Broyden, which is likely due to using estimates for the Jacobian (or its inverse). OPT methods were all able to converge for both cases, albeit slower than the RF methods for the ItP case. Hybrid methods also performed well for both cases, benefiting from the advantages of both classes of methods. \\

In contrast, for a finer mesh consisting of 64 linear elastic elements with a single cohesive element at the same location as the coarse mesh, the performance of each method was related to whether it required an inverse and whether used an estimate for the Jacobian or Hessian. Newton, Dogleg, and both variants of the Broyden method were able to converge for the ItP case.  However, only the variants of the Broyden method were able to converge for the ItC case. This is likely due to the use of estimates for Jacobians, which avoided oscillatory behavior exhibited by methods that use the exact Jacobian. \\

To improve upon the performance of these methods, a promising technique may be restarts (for methods that use a history of iterations) or switching methods depending on recent convergence behavior. Restarts may reintroduce exploratory sub-optimal steps that can escape a local minima or numerically unstable region. Consequently, it can be advantageous to restart a method if convergence is stalling. Similarly, some methods perform well in specific regions of the residual space. Newton (RF) converges rapidly in the neighborhood of the root, while the Broyden methods are more tolerant of local minima. Importantly, as the Krylov methods used for the inner linear solve create smaller and smaller subspaces, there may be eigenvectors and values that the nonlinear solves simply cannot reach; restarting may not help in these situations. \\

For example, a method robust to local minima could be used until iteration reaches the neighborhood of the root at which point the selected method could be switched to Newton, and if Newton exhibits oscillatory behavior, the more robust method could be substituted. Future studies should focus on developing adaptive combinations of methods that leverage the strengths of individual methods and characterize how the methods perform for classical 3D cohesive zone problems, such as a double cantilevered beam. Finally, this study only considers a cohesive zone placed at the location of a node, so the behavior of these methods for problems involving the cohesive segment method and XFEM remains to be characterized.

\section{Acknowledgments}
MKB was supported by the United States Air Force Research Laboratory (AFRL). AC and VS were supported by UES/AFRL subcontract S-200-245-001 to the University of Utah under Prime Contract FA8650-21-D-5279. 

\bibliography{main.bbl}

\end{document}